\documentclass[12pt]{article}
\usepackage{amsmath}
\usepackage{theorem,amsfonts,amssymb}
\newtheorem{theorem}{Theorem}[section]
\newtheorem{proposition}[theorem]{Proposition}
\newtheorem{lemma}[theorem]{Lemma}
\newtheorem{corollary}[theorem]{Corollary}
\theorembodyfont{\upshape}
\newtheorem{definition}{Definition}[theorem]
\newtheorem{example}[theorem]{Example}
\newtheorem{remark}[theorem]{Remark}
\newtheorem{proof}{Proof}

\newtheorem{acknowledgement}{Acknowledgement}

\newcommand{\bt}{\begin{theorem}}
\newcommand{\et}{\end{theorem}}
\newcommand{\bl}{\begin{lemma}}
\newcommand{\el}{\end{lemma}}
\newcommand{\bp}{\begin{proposition}}
\newcommand{\ep}{\end{proposition}}
\newcommand{\bd}{\begin{definition}}
\newcommand{\bex}{\begin{example}}
\newcommand{\eex}{\end{example}}
\newcommand{\ed}{\end{definition}}
\newcommand{\br}{\begin{remark}}
\newcommand{\er}{\end{remark}}
\newcommand{\bc}{\begin{corollary}}
\newcommand{\ec}{\end{corollary}}
\newcommand{\bo}{\begin{proof}}
\newcommand{\eo}{\end{proof}}
\newcommand{\be}{\begin{enumerate}}
\newcommand{\ee}{\end{enumerate}}

\newcommand{\Hom}{{\rm Hom}}

\newcommand{\Aut}{{\rm Aut}}

\newcommand\du{{\rm d}}

\newcommand{\Z}{{\mathbb Z}}
\newcommand{\Q}{{\mathbb Q}}
\newcommand{\N}{{\mathbb N}}

\newcommand{\R}{{\mathbb R}}

\newcommand{\T}{{\mathbb T}}
\newcommand{\Lc}{{\cal L}}

\newcommand{\GL}{{\rm GL}}

\newcommand{\G}{{\cal G}}
\title{Expansive Automorphisms on Locally Compact Groups}
\author{Riddhi Shah}
\date{}

\let\ap=\alpha

\let\ol=\overline

\begin{document}
\maketitle

\begin{abstract}
We show that any connected locally compact group which admits 
an expansive automorphism is nilpotent. 
We also show that for any locally compact group $G$, $\ap\in \Aut(G)$ 
is expansive if and only if for any $\ap$-invariant closed subgroup $H$ which is either compact or normal, 
the restriction of $\ap$ to $H$ is expansive and the quotient map on $G/H$ 
corresponding to $\ap$ is expansive. We get a structure theorem for locally compact groups admitting 
expansive automorphisms. We prove that an automorphism on a non-discrete locally compact group can not be both distal and expansive. 
\end{abstract}


\noindent Secondary: 37B05, 37A25, 22E25





\begin{section}{Introduction} 
Let $G$ be a locally compact group with the identity $e$. An automorphism $\ap$ of $G$ is 
said to be {\it expansive} if $\cap_{n\in\Z} \ap^n(U)=\{e\}$ for some neighbourhood $U$ of $e$; here $U$ is called an 
{\it expansive neighbourhood} for $\ap$. Equivalently, $\ap$ is expansive if there exists a neighbourhood 
$V$ of $e$ such that for every pair $x,y\in G$, $x\ne y$, there exists $n=n(x,y)\in \Z$, such that $\ap^n(y^{-1}x)\not\in V$. 
Expansive automorphisms have been studied extensively on compact groups by 
Kitchens and Schmidt (see \cite{K, KS0, S0} and \cite{S1}), and on totally disconnected groups by Willis \cite{W} and also by Gl\"ockner and Raja \cite{GR}.
The structure of compact groups admitting expansive automorphisms 
has been studied by many (see Kitchens and Schmidt \cite{KS0}, Lam \cite{La} and Lawton \cite{Ln}). 
Any locally compact group admitting an expansive automorphism $\ap$ is metrizable and it has a 
$\sigma$-compact $\ap$-stable open subgroup. Note that if $\ap$ is expansive if and only if $\ap^n$ is so, $n\in \Z$. 

Let $G$ be a compact group and let $\ap\in\Aut(G)$. We say that $(G,\ap)$ satisfies the {\it descending chain condition} (d.c.c.) if for every sequence $G\supset G_1\supset\cdots\supset G_k\supset\ldots$
of closed $\ap$-invariant groups, there exists $K\in\N$ such that $G_k=G_K$ for all $k\geq K$. If $(G,\ap)$ satisfies the d.c.c., then it is conjugate to a Lie subshift, i.e.\ there exist a Lie group $L$ 
and a continuous injective homomorphism $\phi:G\to L^\Z$ such that $\phi\circ\ap=\gamma\circ\phi$, where $\gamma$ is the shift action on $L^\Z$. Moreover, 
$(G,\alpha)$ is said to be a {\it Bernoullian system of Lie type} if  
$\phi$ is bijective. If the shift action in the Bernoullian system $L^\Z$ is expansive then $L$ is finite. An automorphism $\ap$ of a compact group $G$ is expansive, then $(G,\ap)$ satisfies the d.c.c. 
Conversely, if $(G,\ap)$ satisfies the d.c.c.\ and if $G$ is totally disconnected, then $\ap$ is expansive. 

We now describe a solenoidal system (as in \cite{J}). 
Let $U\in \GL(n,\Q)$ and let $\T^m$ denote the $m$-dimensional torus. 
Let $\exp_m:\R^m\to\T^m$ be the exponential function defined as $\exp_m(x)=(\exp 2\pi ix_1,\ldots, \exp 2\pi ix_m)$, for $x=(x_1,\ldots, x_m)\in\R^m$. Let $E_U:\R^m\to(\T^m)^\Z$ be
defined as $E_U(x)=\exp_m(U^n(x))_{n\in\Z}$. The closure of $E_U(\R^m)$ in $(\T^m)^\Z$ is a  shift-invariant subgroup of $(\T^m)^\Z$ which is called a solenoid. Here, $(G,\ap)$ is called a 
{\it solenoidal system} if $G$ is a solenoid and $\ap$ is the shift action on $G$. 

A compact connected group admitting an expansive automorphism is abelian and finite dimensional (see Lam \cite{La} and Lawton \cite{Ln}). 
A locally compact group $G$ is said to be {\it finite dimensional} if it is a 
projective limit of Lie groups, each of which has the same dimension; equivalently if $G$ has a compact totally disconnected normal subgroup $K$ such that $G/K$ is a Lie group. Note that 
$(G,\ap)$ is a solenoidal system if and only if $G$ is a compact connected abelian finite dimensional group and $(G,\ap)$ satisfies the descending chain condition  
 (see \cite{Ln}, \cite{KS0}, \cite{S1} and see also Proposition 6.2 of \cite{J}). In particular, if $G$ is a 
 compact connected group with an expansive automorphism $\ap$, then $(G,\alpha)$ is isomorphic to a solenoidal system, where the matrix $U$ as above does not have any
 of absolute value 1. 

Any compact group $G$ with an expansive automorphism $\ap$ has an open $\ap$-invariant subgroup $H$ of finite index such that $\ap|_H$, the restriction of $\ap$ to $H$, is ergodic (see \cite{S1}). 
 For results on the dynamics of automorphisms on compact groups, we refer the reader to Schmidt  \cite{S1} and the references cited therein. 

Expansive automorphisms on non-compact groups have also been studied (see Eisenberg \cite{Ei}, Aoki \cite{Ao}, Bhattacharya \cite{Bh} and Gl\"ockner and Raja \cite{GR}). 
Theorem 1 in Aoki \cite{Ao} states that any connected locally compact solvable group admitting an expansive automorphism is nilpotent. 
(The present author is however in doubt about whether the proof given in the paper \cite{Ao} is adequate; the argument on page 347, line 6, is unconvincing). 
We show that any connected locally compact group admitting an expansive automorphism is nilpotent (see Theorem \ref{nilp}). 
For a class of compact groups and that of totally disconnected locally compact groups, it is known that expansivity carries over to the quotients modulo closed invariant normal subgroups 
(see Corollary 3.11 in \cite{S0} 
and Theorem A in \cite{GR}). We generalise this to all locally compact groups (see Theorem \ref{quotient}). In addition, we show that the expansivity carries over to the quotient modulo compact invariant 
(not necessarily normal) subgroups (see Theorem \ref{cpt}). We also generalise the structure theorem of Gl\"ockner and Raja (namely, Theorem B of \cite{GR}) for totally disconnected groups to all 
locally compact groups admitting expansive automorphisms (see Theorem \ref{str}). In the end, we show that an automorphism on a locally compact group can not be both distal and expansive 
unless the group is discrete (see Theorem \ref{distal-expa}). 

A homeomorphism $\ap$ of a topological (Hausdorff) space $X$ is said to be {\it distal} if for every pair of distinct elements $x,y\in X$, the closure of $\{(\ap^n(x),\ap^n(y))\mid n\in\Z\}$ in 
$X\times X$ does not intersect the diagonal $\{(g,g)\mid g\in X\}$. If $G$ is a topological group and $\ap\in\Aut(G)$, then the above definition is equivalent to the following: $\ap$ is distal  
if the closure of $\{\ap^n(x)\mid n\in\Z\}$ in $G$ does not contain the identity $e$, for every $x\ne e$. Distal maps on compact spaces were introduced by David Hilbert. Distal automorphisms on locally compact 
groups have been studied extensively. We refer the reader to Raja and Shah \cite{RSh1, RSh2}, Shah \cite{Sh1}, and the references cited therein. Although distality and expansivity are opposite phenomena for 
a non-discrete group, they do satisfy some similar properties. It is easy to see that if the restriction of the automorphism $\ap$ to the closed invariant subgroup and the corresponding map on the quotient have
one of the properties, so does $\alpha$. It has also been shown that distality carries over to the quotients modulo closed invariant subgroups which are either compact or normal 
(cf.\ \cite{RSh1} and \cite {Sh1}).
\end{section}

\begin{section}{Properties and structure of groups with expansive automorphisms}

Any locally compact group $G$ admitting an expansive automorphism $\ap$ has a countable neighbourhood basis of the identity $e$ given by 
$\{\cap_{n=-k}^{k}\ap^n(U)\}_{k\in\N}$, where $U$ is any expansive neighbourhood (of $e$) for $\alpha$ in $G$. Therefore, $G$ is metrizable and it has a left invariant metric 
$d$ compatible with the topology of $G$ (see \cite{HR}). The definitions of expansivity given in the introduction are equivalent to the following definition in terms of the
metric $d$: $\ap$ is expansive if there exists an $\epsilon>0$ such that for every pair $x,y\in G$, $x\ne y$, there exists 
$n=n(x,y)\in\Z$ such that $d(\ap^n(x),\ap^n(y))>\epsilon$. (Note that this is how expansivity for a homeomorphism $\ap$ on a metrizable space $X$ is defined). 

If $H$ is a closed subgroup of $G$, we consider the quotient space $G/H$, the set of left cosets $\{xH\mid x\in G\}$, with the usual quotient topology. 
For $\ap\in \Aut(G)$ and a $\ap$-invariant closed subgroup $H$, we have the canonical map $\bar\ap$ on $G/H$, defined as $\bar\ap(xH)=\ap(x)H$. It is a homeomorphism of
$G/H$. We say that the map $\bar\ap$ on $G/H$ is {\it expansive} if there exists a neighbourhood $U$ of $H$ in $G/H$ such that $\cap_{n\in\Z}\bar\ap^n(U)=H$. 
Equivalently, there exists a neighbourhood $V$ of $H$ in $G/H$ such that for every pair $x,y\in G$ with $xH\ne yH$, there exists $n=n(xH,yH)\in\Z$ such that 
$\bar\ap^n(y^{-1}xH)\not\in V$. This definition coincides with the definition given above for a group when $H$ is normal in $G$ and $G/H$ is a group. 
 If $G$ is metrizable and if a closed subgroup $H$ is normal or compact, then $G/H$ is also metrizable with a $G$-invariant metric; say, $d$, i.e.\ 
 $d(xH, yH)=d(y^{-1}xH,H)=d(gxH, gyH)$ for all $g,x,y\in G$. In this case, the above definition of expansivity on the quotient is equivalent to the 
 one given above in terms of the metric $d$.  It is easy to see that if $\ap$ is an expansive automorphism of $G$ and $H$ is a closed $\ap$-invariant subgroup of $G$, 
 then $\ap|_H$ is expansive. Also, if $\ap|_H$ is expansive and $\bar\ap$ on $G/H$ is expansive, then $\ap$ is expansive on $G$. 

As mentioned above, any compact connected group admitting an expansive automorphism is abelian and finite dimensional (see \cite{La}, \cite{Ln} and also Theorems 5.3 and 6.1 of \cite{KS0}). 
As any connected locally compact group $G$ has a maximal compact normal (characteristic) subgroup $K$ such that $G/K$ is a Lie group, the above also implies that $G$ is finite dimensional if 
it admits an expansive automorphism; (more generally, if $K$ admits an expansive automorphism). For a closed subgroup $H$ of $G$, let $H^0$ denote the connected component of the identity $e$ 
in $H$, it is a closed (normal) characteristic subgroup of $H$. We first state a simple lemma which will be useful.

\begin{lemma} \label{cpt-abelian}
Let $G$ be a connected locally compact group. Let $K$ be any compact normal subgroup of $G$. If $K^0$ is 
abelian, then $K$ is abelian and central in $G$. In particular if $G$, $K$ or $K^0$ admits an expansive automorphism, 
then $K$ is abelian and central in $G$. 
\end{lemma}

\bo
Suppose $K^0$ is abelian. Since $G$ is connected and $K^0$ is abelian and normal in $G$, by Theorem 4 in \cite{Iw},
$K^0$ is central in $G$.  This also implies that $K$ is central in $G$, as given any $g\in G$, there exists 
$h\in K^0$ such that for any $k\in K$, $gkg^{-1}=hkh^{-1}=k$. 

The second assertion follows from the first, as $K^0$ is characteristic in $G$ and $\ap|_{K^0}$ is expansive, we get that
$K^0$ is abelian (cf.\ \cite{La}, Corollary 3.3).
\eo

For a connected Lie group, let $\G$ denote the Lie Algebra of $G$ and let $\exp:\G\to G$ be the exponential map. There is a neighbourhood $U$ of 0 in $\G$ such that $\exp|_U$ is a homeomorphism 
onto the neighbourhood of the identity $e$ in $G$. Let $\du\ap:\G\to\G$ be the Lie algebra isomorphism such that $\exp\circ\,\du\ap(X)=\ap\circ\exp(X)$, $X\in \G$. We note the following which essentially 
follows from Theorem A and Propositions 2.1 and 2.3 of \cite{Bh}. 

\begin{proposition} \label{Lie} Let $G$ be a connected Lie group and let $\ap\in\Aut(G)$. Then the following are equivalent\/{\rm :} 
\begin{enumerate}
\item $\ap$ is expansive.
\item $\du\ap$ is expansive on the Lie algebra $\G$ of $G$. 
\item $\du\ap$ does not have any eigenvalue of absolute 1. 
\end{enumerate}
In particular, if $G$ admits an expansive automorphism, then it is nilpotent.
\end{proposition}

\bo 
(1) $\Rightarrow$ (2) is proven in the proof of Theorem A of \cite{Bh} just by using the fact that one can choose a neighbourhood $U$ of 0 in $\G$ as above such that $\exp(U)$ is an expansive 
neighbourhood for $\ap$. 
(2) $\Leftrightarrow$ (3) follows from Proposition 2.3 of \cite{Bh} (see also \cite{Ei}), and (3) $\Rightarrow$ (1) follows from Proposition 2.1 and Theorem A of \cite{Bh}. 
If $G$ admits an expansive automorphism $\ap$, then $\du\ap$ satisfies condition (3), and $\G$ is a nilpotent Lie algebra 
(see Exercise 21 (b) among the exercises for Part I of \cite{Bo}, \S 4, or Theorem 2 of \cite{Ja}), which in turn implies that $G$ is nilpotent. \eo

Let $G$ be a connected locally compact abelian group and let $\Lc(G)$ denote the space $\Hom(\R,G)$ of all continuous homomorphisms from $\R$ to $G$ 
endowed with the topology of the uniform convergence on the compact subsets of $\R$. Then $\Lc(G)$ is a topological vector space with respect to the pointwise addition and 
scaler multiplication (see Proposition 7.36 in \cite{HM}). If $G$ is a (linear) Lie group, then $\Lc(G)$ coincides with the Lie algebra $\G$ of $G$. In case $G$ is compact, 
$\Lc(G)$ is isomorphic to $\Hom(\hat G,\R)$, where $\hat G$ is the character group of $G$. In general, $G$ is isomorphic to $\R^m\times K$, where $K$ is the maximal (largest) compact 
connected (abelian) subgroup of $G$, and 
$\Lc(G)$ is isomorphic to $\R^m\times \Lc(K)$. Let $\exp:\Lc(G)\to G$ be defined as $\exp(X)=X(1)$, $X\in \Lc(G)$. Then $\exp$ is continuous. For $\ap\in\Aut(G)$, let $\du\ap:\Lc(G)\to\Lc(G)$ be 
defined as $\du\ap(X)=\ap\circ X$, $X\in\Lc(G)$. Note that $\du\ap$ defines a vector space isomorphism of $\Lc(G)$ and $\exp\circ\,\du\ap=\ap\circ\exp$.
Suppose $K$ as above is finite dimensional, then 
$\Lc(G)=\R^m\times\Lc(K)$ is also finite dimensional and it is isomorphic to $\R^{m+n}$, where $\Lc(K)$ is isomorphic to $\R^n$. 
 Moreover, $\exp$ is a continuous homomorphism and the kernel of $\exp$ is contained in $\Lc(K)$, and if $K$ (and hence $G$) is a (linear) Lie group, then $\exp$ is a local isomorphism, i.e.\  
there exists a neighbourhood $U$ of 0 in $\Lc(G)$ (resp.\ $V$ of $e$ in $G$) such that $\exp:U\to V$ is a homeomorphism. If $G$ is isomorphic to $\R^n$
(i.e.\ if $K$ is trivial), then $\exp$ is a vector space isomorphism. 
Note that $\Lc(G)/\Lc(K)$ is isomorphic to $\Lc(G/K)$. As $K$ is connected and $\ap$-invariant, we have $\bar\ap$, the automorphism of $G/K$ corresponding to $\ap$. 
Moreover, $\du\ap$ keeps $\Lc(K)$ invariant and we have that $\ol{\du\ap}$ is the corresponding vector space isomorphism on $\Lc(G)/\Lc(K)$ and $\du\bar\ap=\ol{\du\ap}$ 
(under the isomorphism of $\Lc(G/K)$ and $\Lc(G)/\Lc(K)$). In fact, $G/K$, being a finite dimensional real vector 
space, is isomorphic to $\Lc(G/K)$ under the exponential map. We refer the reader to Ch.\ 7 of \cite{HM} for more details.  
(We use same notations for the exponential map on $\Lc(G)$ and also on the Lie algebra $\G$ of a Lie group $G$. Similarly, we use the same notation for the corresponding vector space isomorphism on 
$\Lc(G)$ as well as for the Lie algebra automorphism when it is induced by an automorphism of a group or a Lie group.) 

As noted earlier, any connected locally compact abelian group $G$ admitting an expansive automorphism is finite dimensional, and hence, so is $\Lc(G)$. Therefore, we can discuss the expansivity of the 
corresponding map on $\Lc(G)$ in the following.  

\begin{lemma} \label{fd} Let $G$ be a connected locally compact abelian group and let $\ap\in \Aut(G)$. Let $K$ be the maximal compact normal subgroup of $G$ and let $\bar\ap$ be the corresponding 
automorphism on $G/K$. If $\ap$ is expansive, then so are $\du\ap$ on $\Lc(G)$ and  $\bar\ap$ on $G/K$. 
\end{lemma}

\bo Suppose $\ap$ is expansive. As observed above, $G=\R^n\times K$, and $K$ as well as $G$ is finite dimensional. Then $\Lc(G)=\R^{m+n}$, where $n=\dim K$ and $\Lc(K)=\R^n$. Let 
$V$ be an expansive neighbourhood of $e$ in $G$ for $\ap$. Let $\exp:\Lc(G)\to G$ be as above. As $\Lc(G)$ is a finite dimensional real vector space, there is a vector space norm on it. As $\exp$ is 
continuous, we can choose $r>0$ such that for the open neighbourhood $U=\{X\in\Lc(G)\mid \|X\|<r\}$ of 0 in $\Lc(G)$, we have $\exp(U)\subset V$. Observe that $tU\subset U$ for all $t\in [-1,1]$.
We show that $U$ is an expansive neighbourhood of $0$ in $\Lc(G)$ for $\du\ap$. 
Let $X\in U$ be such that $\du\ap^n(X)\in U$ for all $n\in \Z$. Then for all $t\in [-1,1]$ and $n\in\Z$, $\du\ap^n(tX)=t\,\du\ap^n(X)\in U$, and hence, $\exp(\du\ap^n(tX))=\ap^n(\exp(tX))\in V$. 
As $V$ is expansive for $\ap$, $\exp(tX)=X(t)=e$ for all $t\in [-1,1]$. Since $X$ is a (real) one parameter subgroup, the preceding assertion implies that $X(t)=e$ for all $t\in\R$, and hence, 
$X=0$. Therefore, $U$ is an expansive neighbourhood of $0$ for $\du\ap$, and hence, $\du\ap$ is expansive on $\Lc(G)$. 

By Proposition 2.3 of \cite{Bh}, the eigenvalues of $\du\ap$ do not have absolute value 1. As $\ap$ keeps $K$ invariant, $\du\ap$ keeps $\Lc(K)$ invariant. 
Let $\ol{\du\ap}$ be the map corresponding to $\du\ap$ on $\Lc(G)/\Lc(K)$. As noted earlier, $\du\bar\ap=\ol{\du\ap}$. Therefore, the eigenvalues of $\ol{\du\ap}$, and hence of $\du\bar\ap$, 
do not have absolute value 1. 
By Proposition 2.3 of \cite{Bh}, $\ol{\du\ap}$, and hence, $\du\bar\ap$ is expansive. Since $G/K$ is a finite dimensional real vector space, the exponential map from $\Lc(G)/\Lc(K)$ to $G/K$ is a vector 
space isomorphism and it is easy to see that the preceding assertion implies that $\bar\ap$ is expansive (see also Proposition \ref{Lie}). \eo 

\begin{remark} It follows from Theorems 8.20 and 8.22 of \cite{HM} that for a connected finite dimensional abelian group $G$, $\exp: \Lc(G)\to G$ is injective on a small neighbourhood of 0, i.e.\ 
given a neighbourhood $V$ of $e$ in $G$, there exists a neighbourhood $U$ of 0 in $\Lc(G)$ such that $\exp(U)\subset V$ and $\exp|_U$ is injective. Hence, for an expansive $\ap\in\Aut(G)$, if $V$ 
is an expansive neighbourhood for $\ap$, then it is easy to see using the injectivity of $\exp|_U$ that $U$ is an expansive neighbourhood for $\du\ap$. This provides an alternative proof for the first 
part of the assertion in Lemma \ref{fd} above.
\end{remark}

\begin{theorem} \label{cpt} Let $G$ be a locally compact group and let $\ap\in\Aut(G)$. Let $K$ be a compact  
$\ap$-invariant subgroup of $G$ and let $\bar\ap:G/K\to G/K$ be the map corresponding 
to $\ap$. If $\ap$ is expansive, then so is $\bar\ap$; equivalently there exists an open set $U$ containing $K$ in $G$ such that 
$\cap_{n\in\Z}\ap^n(U)=K$. 
\end{theorem}

\bo {\bf of Theorem \ref{cpt} for the case when $K$ is normal:} Here $G/K$ is a group and 
$\bar\ap\in\Aut(G/K)$. Suppose 
$K$ is central in $G$. Let $W$ be an expansive neighbourhood of $e$ in $G$ for $\ap$, 
i.e.\ $\cap_{n\in\Z}\ap^n(W)=\{e\}$. Let $V$ be an 
open relatively compact symmetric neighbourhood of the identity $e$ in $G$ such that $\ol{V}^4\subset W$. 
Let $A=\cap_{n\in\Z}\ap^n(V)K=\{x\in V\mid\ap^n(x)\in VK\mbox{ for all }n\in\Z\}K$. Then $K\subset A\subset VK$ and $\ap(A)=A=AK$. If $x,y\in A$, then as 
$K$ is central in $G$, 
$\ap^n(xyx^{-1}y^{-1})\in V^4\subset W$, for all $n\in\Z$ and hence, $xyx^{-1}y^{-1}=e$ since $W$ is an expansive neighbourhood for $\ap$. This implies that 
the elements of $A$, and hence, $\ol{A}$ commute. Let $H$ be the closed subgroup generated by $\ol{A}$. Since $\ol{A}\subset \ol{V}K$ is compact, 
$H$ is compactly generated. Moreover, $H$ is abelian and locally compact. Therefore, $H$ is 
isomorphic to $\R^d\times\Z^k\times C$ and $H^0=\R^d\times C^0$, where $C\subset H$, $C$ is compact and $d,k\in \N\cup\{0\}$. Since $K\subset A$, 
$K\subset C$. Since $\ap(A)=A$, $H$ is $\ap$-invariant and $\ap|_H$ is expansive. Note that $\ap(C)=C$ and $\ap(C^0)=C^0$, as 
$C$ is the maximal compact subgroup of $H$. Since $C$ is compact, the restriction of $\bar\ap$ to 
$C/K$ is expansive (cf.\ \cite{S0}, Corollary 3.11). As
$H^0C/C$ is isomorphic to $H^0/C^0$, by Lemma \ref{fd}, the corresponding action of 
$\ap$ on $H^0C/C$ is expansive. This, together with the preceding assertion, implies that the restriction of $\bar\ap$ to 
$H^0C/K$ is expansive. As $H^0C$ is open in $H$, the restriction of $\bar\ap$ to $H/K$ is also expansive. Let $\pi:G\to G/K$ 
be the natural projection. Let $V'\subset V$ be a neighbourhood of $e$ in $G$ such that $\pi(V')\cap (H/K)$ is an expansive neighbourhood for 
the restriction of $\bar\alpha$ to $H/K$. 
Now we show that $\pi(V')$ is an expansive neighbourhood of $\pi(e)$ for $\bar\ap$  
in $G/K$. Let $x\in V'$ be such that $\pi(x)\in \cap_{n\in\Z}{\bar\ap}^n(\pi(V'))$. Then 
$x\in\cap_{n\in\Z}\ap^n(V)K\subset A$. Therefore, $x\in V'\cap A\subset V'\cap H$. As $\pi(V')\cap (H/K)$ is 
an expansive neighbourhood for the restriction of $\bar\ap$ to $H/K$, we have that $\pi(x)\in\pi(K)$, and hence, $x\in K$. This implies that 
$\pi(V')$ is an expansive neighbourhood of $\pi(e)$ for $\bar\ap$ in $G/K$; i.e.\ $\bar\ap$ is expansive on $G/K$. 

Now suppose $K$ is normal but not central in $G$.  If $G$ is compact, then the assertion follows from Corollary 3.11 of \cite{S0}. 
 Let $K_0$ be the maximal compact normal subgroup of $G^0$. 
Since $\ap|_{G^0}$ is expansive, by Lemma \ref{cpt-abelian}, $K_0$ is central in $G^0$. Here, $K_0$ is a compact characteristic, and hence, $\ap$-invariant normal subgroup in $G$. 
As $G/(KK_0)$ is isomorphic to $(G/K)/((KK_0)/K)$ and the restriction of $\bar\ap$ to $(KK_0)/K$ is expansive 
(cf.\  \cite{S0}, Corollary 3.11), it is enough to prove that the automorphism of $G/(KK_0)$ corresponding to $\ap$ is expansive. Therefore, replacing $K$ by $KK_0$ if necessary, we may
assume that $K_0\subset K$, i.e.\ $G^0\cap K=K_0$.   

Let $W$ be an expansive neighbourhood of the identity $e$ in $G$ for $\ap$, with an additional property that $W$ is contained in an almost 
connected open subgroup of $G$ containing $K$. Let $C$ be a compact subgroup containing $K_0$ such that $C$ is normalised by $G^0$, $CG^0$ is an 
open almost connected subgroup and $(CG^0)/K_0=C/K_0\times G^0/K_0$. 
As $C$ is normal in $CG^0$ which is open and $K$ is compact, $C$ has a subgroup of finite index, which contains $K_0$ and it is normalised by $KG^0$. 
Replacing $C$ by this subgroup, we may assume that $C$ and $K$ normalise each other, and that both $CG^0$ and $CKG^0$ are open almost connected subgroups. Note also 
that $C$ (resp.\ $CK$) is the maximal compact normal subgroup of $CG^0$ (resp.\ $CKG^0$) and $C\cap G^0=K_0=K\cap K_0$. Since $K_0$ is central in $G^0$, 
the automorphism of $G^0/K_0$ corresponding to $\ap$ is expansive. As $(G^0K)/K$ is isomorphic to $G^0/K_0$, the preceding assertion implies that the restriction of $\bar\ap$
to $(G^0K)/K$ is expansive. There exists a symmetric 
relatively compact neighbourhood $U$ of $e$ such that $U\subset CG^0\cap W$ and the image of $U\cap G^0K$ in $(G^0K)/K$ is an expansive 
neighbourhood for the restriction of $\bar\ap$ to $(G^0K)/K$. Let $B=\cap_{n\in\Z}\ap^n(U)K$. Here $K\subset B\subset UK\subset CKG^0$. Suppose $B\ne K$. From our assumption 
on $U$ as above, we get that $B\cap G^0K=K$. Let $H'$ be the closed group generated by $B$ in $G$. Then $H'$ is $\ap$-invariant and  
$K\subset H'\subset CKG^0$; the latter is an open almost connected group. Therefore, $H'G^0$ is also an almost connected group and it is 
$\ap$-invariant. Let $C'$ be the maximal compact normal subgroup of $H'G^0$. Then $C'$ is $\ap$ invariant, $K\subset C'$ and $H'G^0\cap CK\subset C'$.   
Therefore, $H'G^0=(H'G^0\cap CK)G^0\subset C'G^0$. This implies that $H'G^0=C'G^0$ and  
 $(H'G^0)/K=C'/K\times (G^0K)/K$. Note that $C'$ is $\ap$-invariant and $\ap|_{C'}$ is expansive. As $C'$ is compact, 
by corollary 3.11 of \cite{S0}, the restriction of $\bar\ap$ on $C'/K$ is expansive. As the restriction of $\bar\ap$ to $(G^0K)/K$ is also
expansive, the preceding assertion implies that the restriction of $\bar\ap$ to $(H'G^0)/K$ is expansive. In particular, the restriction of $\bar\ap$ to $H'/K$ is expansive. 

Let $U'\subset U$ be a neighbourhood of the identity $e$ in $G$ such that the image of $U'\cap H'$ in $H'/K$ is an expansive 
neighbourhood for the restriction of $\bar\ap$ to $H'/K$. Arguing as in the first part of the proof by replacing $V'$, $A$ and $H$ by $U'$, $B$ and 
$H'$ respectively, it is easy to deduce that the image of $U'$ in $G/K$ is an expansive neighbourhood for $\bar\ap$. This completes the proof for the case when $K$ is normal. \eo

We will complete the proof of Theorem \ref{cpt} after the next result. There are examples of expansive automorphisms on connected nilpotent (Lie) groups; see an example in \S 3 of \cite{Bh}. 
One can also produce a connected nilpotent group which is not  a Lie group by taking a product of a connected nilpotent Lie group with an expansive automorphism and a particular solenoid 
(which is not a Lie group) whose shift map is expansive. The following theorem shows that there is no connected locally compact non-nilpotent group which admits expansive automorphisms. 

\begin{theorem} \label{nilp} Any connected locally compact group admitting an expansive automorphism is nilpotent. 
\end{theorem}

\bo Let $G$ be a connected locally compact group and let $\ap\in \Aut(G)$ be expansive. Let $K$ be the maximal compact normal subgroup of $G$. By Lemma \ref{cpt-abelian}, 
$K$ is abelian and central in $G$. Moreover, we get from Theorem \ref{cpt} for the normal case (proven above) that the automorphism on $G/K$ corresponding to $\ap$ is expansive. 
Note that $G$ is nilpotent if $G/K$ is so. Therefore, it is enough if we assume that $G$ is a connected Lie group without any nontrivial compact normal subgroup. 
Now the assertion follows from Proposition \ref{Lie}. \eo

\bo {\bf of Theorem \ref{cpt} for the general case:} Here, the compact group $K$ is not assumed to be normal in $G$. By Theorem 2.3,  $G^0$ is nilpotent. 
Then any compact subgroup of $G^0$ is abelian and central in $G^0$. 
Let $K_0$ be the maximal compact normal subgroup of $G^0$. Then $K_0$ 
is characteristic in $G$, and hence it is $\ap$-invariant and normal in $G$. Note that as $K_0$ is abelian and normal in $G$, 
$K\cap K_0$ is normal in $KK_0$. By Corollary 3.11 of \cite{S0}, the $\ap$-action on $(KK_0)/(K\cap K_0)$ is expansive. Let $\pi_1:KK_0\to (KK_0)/(K\cap K_0)$ be the natural projection. Then 
$\pi_1(KK_0)=\pi_1(K)\ltimes \pi_1(K_0)$, a semidirect product. One can choose an expansive neighbourhood $V$ in $\pi_1(K_0)$ for the $\ap$-action on 
$\pi_1(K_0)$, and it follows that the image of $V\pi_1(K)$ in $\pi_1(KK_0)/\pi_1(K)$ is an expansive neighbourhood of $\pi_1(K)$ for the $\ap$-action on $\pi_1(KK_0)/\pi_1(K)$. Therefore, the 
$\ap$-action on $\pi_1(KK_0)/\pi_1(K)$ is expansive. Since the latter is isomorphic to 
$(KK_0)/K$, we have that the $\ap$-action is expansive on $(KK_0)/K$. Now it is enough to prove that the $\ap$-action is expansive on $G/(KK_0)$. 
Note that $G/(KK_0)$ is isomorphic to $(G/K_0)/((KK_0)/K_0)$ and, as $K_0$ is normal, from the proof of the normal case above, the $\ap$-action on $G/K_0$ is 
expansive. Therefore, replacing $G$ by $G/K_0$ and $KK_0$ by $(KK_0)/K_0$, without loss of any 
generality, we may assume that $K_0$ is trivial. Now $G^0$ is a simply connected nilpotent Lie group without any non-trivial compact subgroup. In particular, 
$K$ is totally disconnected. 

Let $U$ be a relatively compact expansive neighbourhood of $e$ in $G$ for $\ap$. We first show that there exists a sequence of compact totally disconnected groups $\{C_n\}_{n\in\N}$ such that 
$\cap_n C_n=\{e\}$, and for each $n$, $C_{n+1}\subset C_n\subset U$, $C_nG^0$ is an open subgroup in $G$ and each $C_n$ is normalised by $K$ as well as $G^0$. 
As $\ap$ is expansive, we have that $G$ is metrizable, and hence, so is $G/G^0$.  Since $K$ is compact and $G/G^0$ is totally disconnected, 
 $G/G^0$ has a neighbourhood basis $\{B_n\}_{n\in\N}$ consisting of compact open subgroups such that $B_{n+1}\subset B_n$ and each $B_n$ is invariant under the 
 conjugation action of $K$ on $G/G^0$. 
 Since each $B_n$ has finite index in $B_1$, $B'_n= \cap_{b\in B_1}bB_nb^{-1}$ is an open subgroup of finite index in $B_n$ and it is also invariant
 under the conjugation action of $K$ on $G/G^0$. Replacing $B_n$ by $B'_n$ if necessary, we may assume that each $B_n$ is normal in 
 $B_1$. Therefore, $G$ itself has open subgroups 
$H_n$ such that $H_n/G^0=B_n$, $H_{n+1}\subset H_n$, for all $k\in K$, $kH_nk^{-1}=H_n$ for all $n$ and $\cap_n H_n=G^0$. Moreover, each $H_n$ is normal in $H_1$. 
As each $H_n$ is almost connected, it has a unique maximal compact normal subgroup $C_n$ such that $K$ normalises $C_n$ and $H_n/C_n$ is a Lie group with finitely 
many connected components. In particular, each $C_nG^0$ is an open subgroup in $G$ and $C_n\cap G^0=\{e\}$ as $G^0$ has no nontrivial compact subgroup. 
Here, $G^0$ normalises $C_n$, since $G^0\subset H_n$. Therefore, $C_nG^0=C_n\times G^0$ and, $C_n$ and $C_nG^0$ are normalised by $K$.  
Since $\cap_n C_n\subset \cap_n H_n=G^0$, $\cap_n C_n=\{e\}$. Observe that each $C_n$ is 
characteristic in $H_n$, and hence, it is normal in $H_1$. Therefore, $C_{n+1}\subset C_n$ for all $n$. This, together with the preceding assertion, impliies that 
$C_n\subset U$ for all large $n$, (say) for al $n\geq n_0$. Replacing $C_n$ by $C_{n+n_0}$ for each $n$, we get that $C_n\subset U$ for all $n$.
Since $C_n\cap G^0=\{e\}$, each $C_n$ is totally disconnected. 

Let $U_n=C_n\times W_n\subset U$, $n\in\N$. Here, $\{U_n\}_{n\in\N}$ is a 
neighbourhood basis of the identity $e$ in $G$ such that $U_{n+1}\subset U_n$ and $W_n=U_n\cap G^0$ is open in $G^0$, $n\in\N$. Fix any $n\in \N$. Suppose 
$x\in U_n$ is such that $\ap^m(x)\in U_nK$, for all $m\in\Z$. Then $x=wc=cw$ for some $c\in C_n$ and $w\in W_n$. Now for $m\in\Z$, 
$\ap^m(x)=c_mw_mk_m=\ap^m(c)\ap^m(w)=\ap^m(w)\ap^m(c)$, where 
$c_m\in C_n$, $w_m\in W_n$, $k_m\in K$, $c_0=c$, $w_0=w$ and $k_0=e$.  As both $C_n$ and $\ap^m(C_n)$ centralise $G^0$, we get that 
$$\ap^m(c^{-1})c_mk_m=w_m^{-1}\ap^m(w)=c_mk_m\ap^m(c^{-1}).$$ 
As $C_nK$ and $\ap^m(C_n)$ are compact subgroups and $\ap^m(c^{-1})$ and $c_mk_m$ 
commute with each other, we get that $w_m^{-1}\ap^m(w)=\ap(c^{-1})c_mk_m$ generates a compact subgroup in $G^0\cap\ap^m(C_n)C_nK_n$. As $G^0$ has no 
nontrivial compact subgroup, we get that $w_m^{-1}\ap^m(w)=e$, and hence that $\ap^m(w)=w_m\in W\subset U$. Since this holds for all $m\in\Z$ and since $U$ is expansive for $\ap$, 
$w=e=w_m$ and $x\in C_n$. Now $\ap^m(x)=c_mk_m$, for all $m\in\Z$. 
Let $C'_n=\{c\in C_nK\mid\ap^m(c)\in C_nK\mbox{ for all }m\in\Z\}$. Then $x\in C'_n$. For each $n$, $C'_n$ is a closed (compact) subgroup of $C_nK$, $K\subset C'_{n+1}\subset C'_n\subset C'_1$ 
and $\cap_nC'_n=K$. Each $C'_n$ is a $\ap$-invariant compact subgroup of $C'_1$. As $\ap|_{C'_1}$ is  expansive, it satisfies the descending chain condition (d.c.c.), and hence, 
there exists $N\in\N$ such that for all $n\geq N$, $C'_N=\cap_n C'_n=K$. The above implies that if $x\in U_N=C_N\times W_N$  
is such that $\ap^m(x)\in U_NK$ for all $m\in\Z$, then $x\in C'_N=K$. This shows that $\bar\ap$ on $G/K$ is expansive. 
\eo

The following theorem generalises Theorem A of \cite{GR} for totally disconnected locally compact group to all locally compact groups. 

\begin{theorem} \label{quotient}
Let $G$ be a locally compact group and let $\ap\in\Aut(G)$. Let $H$ be a closed normal $\ap$-invariant subgroup of $G$ and let 
$\bar\ap$ be the automorphism of $G/H$ corresponding to $\ap$. Then $\ap$ is expansive if and only if $\ap|_H$ and
$\bar\ap$ are expansive. 
\end{theorem}

\bo
The `if' statement follows easily. Now suppose $\ap$ is expansive. Let $H$ be a closed normal $\ap$-invariant subgroup of $G$. 
Clearly $\ap|_H$ is expansive. We show that $\bar\ap$ on $G/H$ is expansive. 
If $H$ is compact or $G$ is totally disconnected, then the assertion 
follows from Theorem \ref{cpt} above or Theorem A of \cite{GR} respectively. Let $K$ be the maximal compact normal subgroup of $G^0$.  Then $HK$ is a 
closed normal $\ap$-invariant subgroup and $(HK)/H$ is isomorphic to $K/(K\cap H)$. Since the $\ap$-action on $K$ is expansive, 
so is the corresponding action on $K/(K\cap H)$ (cf.\  \cite{S0}, Corollary 3.11). Hence, it is enough to show that the $\ap$-action on $G/HK$ is expansive, i.e.\ we may assume that 
$K\subset H$. Moreover, $G/H$ is isomorphic to $(G/K)/(H/K)$ and from Theorem \ref{cpt}, the $\ap$-action on $G/K$ is expansive. Replacing $G$ by 
$G/K$ and $H$ by $H/K$, we may assume that $G^0$ has no nontrivial compact  normal subgroup. As $\ap$ is expansive, by  Theorem \ref{nilp}, $G^0$ is nilpotent, and hence, a 
simply connected nilpotent Lie group.
Suppose $H$ is connected. Then $H\subset G^0$ and both $H$ and $G/H$ are also simply connected and nilpotent. Let $\G$ (resp.\ $\mathcal H$) denote the 
Lie algebra of $G^0$ (resp.\ $H$) and let $\ap_0=\ap|_{G^0}$. By Proposition \ref{Lie}, the eigenvalues of $\du\ap_0$ on $\G$ do not have absolute value 1. The same holds 
for the eigenvalues of $\du\bar\ap_0$, the corresponding map on the Lie algebra $\G/\mathcal H$ of $G^0/H$. By Proposition \ref{Lie}, $\bar\ap_0$, 
the restriction of $\bar\ap$ to $G^0/H$ is expansive. 

Let $U$ be an expansive relatively compact symmetric neighbourhood of the identity $e$ in $G$ for $\ap$. As $G^0$ has no nontrivial compact normal subgroup, there exists a 
compact totally disconnected subgroup 
$C$ normalised by $G^0$ such that $C\subset U$, $C\cap G^0=\{e\}$ and $C\times G^0$ is open in $G$. Replacing $U$ by a smaller symmetric neighbourhood, 
we may assume that $U=C\times W$, where $W=U\cap G^0=W^{-1}$, and that the image of $W$ in $G^0/H$ is an expansive neighbourhood for $\bar\ap_0$. 

Let $x\in U$ be such that $\ap^n(x)\in UH$ for all $n\in\Z$.  As $UH=CWH\subset CG^0$, we have $x=cw$ and 
$\ap^n(c)\ap^n(w)=\ap^n(x)=c_nw_nh_n\in C\times G^0$, where $c,c_n\in C$, $w,w_n\in W$ and $h_n\in H$ for all $n\in\Z$, $c_0=c$, $w_0=w$ and $h_0=e$. Fix any $n\in\N$. 
Using the fact that both $C$ and $\ap^n(C)$ centralise $G^0$, from above we get that $\ap^n(c)$ and $c_n$ commute with each other. Now the above implies that 
$\ap^n(c)c_n^{-1}=\ap^n(w^{-1})w_nh_n$ generates a compact subgroup in $\ap^n(C)C\cap G^0$. Since $G^0$ has no nontrivial 
compact subgroup, we get that $\ap^n(c)=c_n$ for all $n\in\Z$. As $U$ is expansive for $\ap$, we have that $c=e=c_n$, and hence, that $x=w$. Therefore, $\ap^n(w)=w_nh_n$, $n\in\Z$. As the image 
of $W$ in $G^0/H$ is an expansive neighbourhood for $\bar\ap_0$, the restriction of $\bar\ap$ to $G^0/H$, the preceding assertion implies that $w\in H$. This proves that $\bar\ap$ is expansive 
if $H$ is connected. 

Now suppose $H$ is not connected. Observe that $H^0$ is $\ap$-invariant and normal in $G$. From above, we have that the $\ap$-action on $G/H^0$ is expansive. Since 
$G/H$ is isomorphic to $(G/H^0)/(H/H^0)$, we may replace $G$ by $G/H^0$ and $H$ by $H/H^0$ and assume that $H^0=\{e\}$ 
and that $H$ is totally disconnected. Here, $G^0$ has no nontrivial compact subgroup and it is a simply connected nilpotent subgroup. 

Let $C$ and $U=C\times W$ be as above. Observe that $CH$ is a closed subgroup. 
As $C$ is a compact totally disconnected group, it has arbitrarily small compact open normal subgroups $C_n$ such that $C_nG^0=C_n\times G^0$ is open in $G$,
$C_nG^0$ (resp.\ $C_nH$) are open normal subgroups in $CG^0=C\times G^0$ (resp.\ $CH$) and $\cap_n C_n=\{e\}$.  This, together with the fact that $H$ is totally disconnected, implies that
$(CH)^0=H^0=\{e\}$, and hence, that $CH$ is also totally disconnected. Let $C_H=\{c\in CH\mid \ap^n(c)\in CH\}$. Then $C_H$ is a closed $\ap$-invariant group 
and $H\subset C_H$ is co-compact. As $CH$ is totally disconnected, so is $C_H$. By Theorem A in \cite{GR}, the restriction of $\bar\ap$ to $C_H/H$ is expansive. Let $U'\subset U$ be a 
neighbourhood of the identity $e$ such that the image of $U'H\cap C_H$ in $C_H/H$ is an expansive 
neighbourhood for the restriction of $\bar\ap$ to $C_H/H$. 

Since $\ap$ is continuous, we can choose an open symmetric neighbourhood $V$ of the identity $e$ in $G$ of the form $V=C'\times W'$, $C'=C_{n_0}$ for some fixed $n_0\in\N$, such that  
$V\ap(V)\ap^{-1}(V)\subset  U'\subset U$, 
where $W'$ is an open symmetric (relatively compact) neighbourhood of $e$ in $G^0$. As $CH$ is totally disconnected, we can choose 
$V$ with an additional property that $V\ap(V)\ap^{-1}(V)\cap CH$ is contained in an  
open compact subgroup of $CH$. Since $G^0$ has no nontrivial compact subgroup, we get that $W'\ap(W')\ap^{-1}(W')\cap CH=\{e\}$. As
$C'\ap(C')\ap^{-1}(C')\subset U=C\times W$, for any $x\in C'$, $\ap(x)=cw=wc$ for  some $c\in C$ and $w\in W\subset G^0$. 
This implies that $c^{-1}\ap(x)=\ap(x)c^{-1}$ generates a compact group in 
$C\ap(C')\cap G^0$. As $G^0$ has no nontrivial compact subgroup, $c^{-1}\ap(x)=e$. This shows that $\ap(C')\subset C$ and $C'\ap(C')\subset C$. 
Similarly, we get that $C'\ap^{-1}(C')\subset C$. 

Let $x\in V$ be such that $\ap^m(x)\in VH$, $m\in\Z$. Then $x=cw$ and 
$\ap^m(x)=\ap^m(c)\ap^m(w)=c_mw_mh_m$, where $c,c_m\in C'$, 
$w,w_m\in W'$ and $h_m\in H$, $m\in\Z$, $c_0=c$, $w_0=w$ and $h_0=e$. We show that $\ap^m(w)=w_m$ and $\ap^m(c)\in c_mH$ for all $m\in\N$ by induction. For $m=1$,
$\ap(c)\ap(w)=c_1w_1h_1$, and hence, $c_1^{-1}\ap(c)=w_1h_1\ap(w^{-1})=w_1\ap(w^{-1})h'_1$ for some $h'_1\in H$; which exists as $H$ is normal. Therefore,   
$w_1\ap(w^{-1})\in W'\ap(W')\cap C'\ap(C')H\subset W'\ap(W')\cap CH=\{e\}$, from the choice of $V$ and $W'$ as above; and hence,
 $\ap(w)=w_1$. Now $\ap(c)=c_1h'_1$, i.e.\ $\ap(c)\in c_1H$. 
For a fixed $k\in\N$, suppose $\ap^k(w)=w_k$ and $\ap^k(c)\in c_kH$. We have 
$$
\ap^{k+1}(x)=\ap(\ap^k(c)\ap^k(w))=\ap(c_k)\ap(w_k)\ap(h_k)=c_{k+1}w_{k+1}h_{k+1}.$$ 
This implies that $\ap(c_k)\ap(w_k)=c_{k+1}w_{k+1}h'$ for some $h'\in H$. 
Arguing as above for $c_k,w_k$ instead of $c,w$ and $c_{k+1}, w_{k+1}, h'$ instead of $c_1, w_1, h_1$ respectively, we get that 
$\ap(w_k)=w_{k+1}$ and $\ap(c_k)\in c_{k+1}H$. Using the induction hypothesis for $k$, we have
$\ap^{k+1}(w)=\ap(w_k)=w_{k+1}$ and $\ap^{k+1}(c)\in\ap(c_k)H= c_{k+1}H$. This proves the statement for all $m\in\N$ by induction. Replacing $\ap$ by $\ap^{-1}$ and using the 
facts that $C'\ap^{-1}(C')\subset C$ and $W'\ap^{-1}(W')\cap CH=\{e\}$, we get as above, that $\ap^{-m}(w)=w_{-m}$ and $\ap^{-m}(c)\in c_{-m}H$ for all 
$m\in\N$. Since $\ap^m(w)\in W'\subset U$, $m\in\Z$, and $U$ is expansive for $\ap$, it follows that $w=e$, and hence, $x\in C'$ and $\ap^m(x)\in C'H$, for all $m\in\Z$. Since $C'\subset C$, 
we have that $x\in C'\cap C_H$ and 
$\ap^n(x)\in C'H\cap C_H$ for all $n\in\Z$. As $C'\subset V\subset U'$, and the image of $U'H\cap C_H$ in $C_H/H$ is expansive for the restriction of $\bar\ap$ to 
$C_H/H$, we get that $x\in H$. This implies that the image of $VH$ in $G/H$ is 
an expansive neighbourhood for $\bar\ap$. Therefore, $\bar\ap$ is expansive.
\eo

Note that  Theorem \ref{cpt} implies that Theorem \ref{quotient} holds also when $H$ is a compact (not necessarily normal) subgroup. 

\begin{remark} It is well-known that any compact connected group admitting an expansive automorphism is abelian, finite dimensional and isomorphic to a solenoid, i.e.\ if $G$ is compact 
and connected and $\ap\in\Aut(G)$ is expansive then $(G,\ap)$ is a solenoidal system (cf.\ \cite{Ln}, \cite{KS0}, \cite{S1}). 
The above also implies that any almost connected  locally compact group $G$ admitting an expansive automorphism is finite dimensional. 
\end{remark}

A locally compact group is said to be {\it topologically perfect} if  its commutator subgroup is dense in the whole group. Theorem B of \cite{GR} describes the structure of a 
totally disconnected locally compact group admitting an expansive automorphism. We generalise this to all locally compact groups. 

\begin{theorem} \label{str} Let $G$ be a locally compact group and $\ap\in Aut(G)$ be expansive. Then there exist $\ap$-invariant closed subgroups 
$$G=G_0\supseteq G_1\supseteq\cdots \supseteq G_n=\{e\}$$ 
of $G$ such that $G_j$ is normal in $G_{j-1}$ for $j \in\{1, . . . , n\}$ and each of the quotient groups $G_{j-1}/G_j$ is discrete, abelian or topologically perfect. 
Moreover,  one can choose $\{G_j\}$ in such a way that every $\ap_j$-invariant closed normal subgroup of 
$G_{j-1}/G_j$ is discrete or open, where $\ap_j : G_{j-1}/G_j\to G_{j-1}/G_j$ is defined as $gG_j \mapsto\ap(g)G_j$ for all $g\in G_{j-1}$, for all $j$.
\end{theorem}

\bo Let $\bar\ap:G/G^0\to G/G^0$ be the automorphism of $G/G^0$ corresponding to $\ap$. Then 
by Theorem \ref{quotient}, $\bar\ap$ is expansive. As $G/G^0$ is totally disconnected, by Theorem B of \cite{GR}, the assertion holds for $G/G^0$ and $\bar\ap$. Now it is enough to 
show that the assertions holds for a connected group $G$. We produce a finite sequence of closed normal subgroups satisfying the first assertion, which would in turn imply the first 
assertion for a general $G$. Here, $G$ is a (connected) nilpotent group (cf.\ Theorem \ref{nilp}). Let $K$ be the maximal compact normal subgroup of $G$, then $K$ is characteristic 
in $G$ and it is connected and central, as $G$ is connected and nilpotent. Moreover, $G/K$ is a simply connected nilpotent Lie group. Observe that $G$ has a central series of connected 
characteristic subgroups $G^{(1)}=\ol{[G,G]}$, the closure of the commutator subgroup of $G$, and $G^{(m+1)}=\ol{[G, G^{(m)}]}$ such that for some $k\in\N$, $G^{(k-1)}$ is nontrivial 
and $G^{(k)}$ is trivial. Let $G_0=G$, $G_m=G^{(m)}K$, $m\in\{1,\ldots, k\}$ and let 
$G_{k+1}=\{e\}$. We have that $G$ has a finite sequence of closed normal decreasing subgroups $\{G_m\}$ whose successive quotients, being simply connected and abelian, 
are isomorphic to $\R^{n_m}$ (for some $n_m\in\N$) except for the last one, which is isomorphic to $K$, which is abelian. Therefore the first assertion in the theorem holds. 
Now we expand this finite sequence to a possibly larger finite sequence of closed normal subgroups for which the second assertion 
in the theorem also holds. As $G_{m-1}/G_m$ is central in $G/G_m$ for all $m$, we have that any subgroup of $G_{m-1}$, which contains $G_m$, is normal in the connected group $G$. 
Since the automorphism corresponding to $\ap$ on each $G_{m-1}/G_m$ is expansive (cf.\ Theorem \ref{quotient}), it is enough to assume that $G=\R^n$ or $G=K$, a 
compact connected abelian group. 

Now suppose $G=\R^n$ and $\ap$ is an invertible linear map. We take a sequence of $\ap$-invariant subspaces $\{V_j\}$ such that $V_0=G$ and for $j\geq 1$, if $V_{j-1}\ne\{e\}$, then 
$V_j$ is the maximal closed proper $\ap$-invariant subspace in $V_{j-1}$. This implies that $\dim V_{j-1}>\dim V_j$, Therefore, there exists $k$ such that $V_k=\{e\}$, and hence, the 
sequence $\{V_j\}$ is finite. If $\ap_j:V_{j-1}/V_j\to V_{j-1}/V_j$ is the natural quotient map defined from the restriction of $\ap$ to $V_{j-1}$, then from the choice of $\{V_j\}$ as above, 
any $\ap_j$-invariant subgroup is either discrete or whole of $V_{j-1}/V_j$. So far the expansivity in the connected group is used only to ascertain that it is nilpotent, the assertions in the 
theorem would follow for any (not necessarily expansive) automorphism for a simply connected nilpotent group. 

Now suppose $G=K$, a compact connected abelian group. Then $G$ is finite dimensional. Let 
$G_0=G$ and for $j\geq1$, if $G_{j-1}\ne\{e\}$, then we choose $G_j$ to be the largest proper closed connected (compact) $\ap$-invariant subgroup of $G_{j-1}$. As $G$ is 
finite dimensional and abelian, it is possible to choose such a sequence of $\{G_j\}$. Moreover, as $(G,\ap)$ satisfies the d.c.c.\ (or as $G$ is finite dimensional), we have 
that there exists $k$ such that $G_k=\{e\}$. For $\ap_j$ defined as above, $\ap_j$ is expansive and $(G_{j-1}/G_j,\ap_j)$ is a solenoidal system. Due to the choice of $G_j$, 
we have that any proper closed $\ap_j$-invariant subgroup in $G_{j-1}/G_j$ is totally disconnected,
and hence, it is finite (cf.\ \cite{J}, Theorems 6.4). This completes the proof. \eo

Recall that an automorphism $\ap$ of a locally compact group $G$ is distal if the closure of $\{\ap^n(x)\mid n\in\Z\}$ does not contain the identity $e$ for every $x\ne e$. 

\begin{theorem} \label{distal-expa}
Let $G$ be a locally compact group and let $\ap\in\Aut(G)$. Then $G$ is discrete if and only if $\ap$ is both expansive and distal.
\end{theorem}

\bo The `only if' statement is obvious. Now suppose $\ap$ is both expansive and distal.
 Suppose $G$ is a connected Lie group.  As in the proof of Theorem \ref{nilp}, $G$ is nilpotent. As $\ap$ is expansive, either $G$ is trivial or 
 $\du\ap$, the corresponding Lie algebra automorphism on the Lie algebra of $G$, do not have any eigenvalue of absolute value 1 (cf.\ Proposition \ref{Lie}). On the other hand, since $\ap$ is distal, 
 all the eigenvalues of $\du\ap$ have absolute value 1 (cf.\ \cite{Ab1, Ab2}). This implies that $G$ is trivial. Suppose $G$ is a compact group.  Since $\ap$ is expansive, $(G,\ap)$ satisfies the 
 descending chain condition. As $\ap$ is distal, the preceding assertion implies that $G$ is a Lie group (cf.\ \cite{RSh2}, Lemma 2.4). 

Let $K$ be the maximal compact normal subgroup of $G^0$. Then $K$ is characteristic in $G$, $G^0/K$ is a Lie group 
and $\ap|_K$ is expansive as well as distal. From above, $K$ is Lie group. Hence $G^0$ itself is a Lie group. As $\ap|_{G^0}$ is expansive as well as distal, we get that $G^0$ is trivial, and hence, 
$G$ is totally disconnected. 
As $\ap$ is distal, by Proposition 2.1 of \cite{JR},  $G$ has a neighbourhood basis of compact open $\ap$-invariant subgroups. As $\ap$ is expansive, it leads to a contradiction unless $G$ is discrete. 
\eo
\end{section}
\begin{acknowledgement} The author would like to thank Helge Gl\"ockner for an extensive discussion and the Fields 
Institute, Toronto, Canada for hospitality while some part of the work was done as a visiting scientific researcher to the 
{\it Theme Period on Group Structure, Group Actions and Ergodic Theory} in February 2014. 
\end{acknowledgement}

\bigskip\smallskip
\advance\baselineskip by 2pt
\noindent Riddhi Shah \\
School of Physical Sciences(SPS)\\
Jawaharlal Nehru University(JNU)\\
New Delhi 110 067, India\\
rshah@mail.jnu.ac.in\\
riddhi.kausti@gmail.com

\end{document}